\documentclass[12pt]{amsart} 
\title[Hyperbolic groups with carpet boundaries]{
Local rigidity for hyperbolic groups with Sierpi\'nski carpet boundaries
}
\author{Sergei Merenkov}
\address{Department of Mathematics\\
University of Illinois\\ 1409 W Green St\\ Urbana, IL
61801\\USA} \email{merenkov@illinois.edu}
\thanks{Supported by NSF grant DMS-1001144}

\subjclass[2010]{52C25}
\date{\today}

\usepackage[active]{srcltx}

\usepackage{amsmath, amssymb, amsthm,latexsym}
\usepackage{graphicx}
\newcommand\C{{\mathbb C}}

\newcommand\N{{\mathbb N}}

\renewcommand\H{{\mathbb H}}

\newcommand\Ju{{\mathcal J}}
\newcommand\dee{\partial}

\renewcommand\:{\colon}

\newcommand\hC{{\hat{\mathbb C}}}
\newcommand\no{\noindent}
\newcommand\crit{\rm{crit}}

\newtheorem{theorem}{Theorem}[section]

\newtheorem{corollary}[theorem]{Corollary}

\newtheorem{lemma}[theorem]{Lemma}

\theoremstyle{definition}

\begin{document}


\abstract
{
Let $G$ and $\tilde G$ be Kleinian groups  whose limit sets  $S$ and $\tilde S$, respectively, are homeomorphic to the standard Sierpi\'nski carpet, and such that every complementary component of each of $S$ and $\tilde S$ is a round disc. 
%
We assume that the groups $G$ and $\tilde G$ act cocompactly on triples on their respective limit sets. The main theorem of the paper states that any quasiregular map (in a suitably defined sense) from an open connected subset of $S$ to $\tilde S$  is the restriction of a M\"obius transformation that takes $S$ onto $\tilde S$, in particular it has no branching. This theorem applies to the fundamental groups of compact hyperbolic 3-manifolds with non-empty totally geodesic boundaries. 

One consequence of the main theorem is the following result. Assume that $G$ is a torsion-free hyperbolic group whose boundary at infinity $\dee_\infty G$ is a Sierpi\'nski carpet that embeds quasisymmetrically into the standard 2-sphere. Then there exists a group $H$ that contains $G$ as a finite index subgroup and such that any quasisymmetric map $f$ between open connected subsets of $\dee_\infty G$ is the restriction of the induced boundary map of an element $h\in H$.   
}
\endabstract

\maketitle

\section{Preliminaries}\label{s:Intro}
%
%

\no
Let $\hC$ denote the Riemann sphere. Whenever needed, we identify $\hC$ with the ideal boundary of hyperbolic 3-space $\H^3$. We also assume that $\hC$ is equipped with the chordal metric and call this metric space the \emph{standard 2-sphere}.

\subsection{Group actions}
Recall that a \emph{Kleinian group} $G$ is any discrete group of orientation preserving isometries of $\H^3$. Equivalently, it is a discrete group of orientation preserving M\"obius transformations of the Riemann sphere $\hC$. The \emph{limit set} $\Lambda(G)\subseteq\hC$ of a Kleinian group $G$ is the set of accumulation points of the orbit $G_p$ of any element $p$ in $\H^3$. The limit set $\Lambda(G)$ of any infinite Kleinian group $G$ is a non-empty compact subset of $\hC$ and the group $G$ acts on $\Lambda(G)$ by homeomorphisms. If the Kleinian group $G$ is non-elementary, i.e., its limit set has more than two points, then $\Lambda(G)$ is a  perfect set. 

Let $Y$ be a locally compact Hausdorff topological space. An action of a discrete group $G$ on $Y$ by homeomorphisms is said to be \emph{properly discontinuous} 
if for all compact subsets $K$ and $L$ of $Y$ the set
$$
\{g\in G\: g(K)\cap L\neq\emptyset\}
$$
is finite. Such an action is called \emph{cocompact} if the quotient $Y/G$ is compact.

If $X$ is any compact Hausdorff topological space that has at least three points, we denote the space of distinct triples of $X$ by $\Sigma_3(X)$, namely
$$
\Sigma_3(X)=\{(o,p,q)\in X^3\: o\neq p, o\neq q, p\neq q\}.
$$ 
Assume that a group $G$ acts on $X$ by homeomorphisms.
Such an action induces a  diagonal action of $G$ on $\Sigma_3(X)$.
Following~\cite{BB96}, we say that a group $G$ is a \emph{uniform convergence group} acting on a perfect, compact, Hausdorff topological space $X$ if the action of $G$ on $\Sigma_3(X)$ is properly discontinuous and cocompact. A cocompact action of a discrete group on the space of distinct triples is referred to as \emph{cocompact on triples}.

\subsection{Hyperbolic groups}\label{SS:HG}
A finitely generated group $G$ is called \emph{hyperbolic} if there exists a symmetric finite generating set $S$ for $G$ and a positive constant $\delta$ such that the geodesic triangles of the Cayley graph of $G$ with respect to $S$ are $\delta$-\emph{thin}. The latter means that any side of a geodesic  triangle is contained in the $\delta$-neighborhood (with respect to the word metric) of the union of the other two sides. See, e.g.,~\cite{gh}, for  background on hyperbolic groups. Also, see~\cite{KB01} for a survey on hyperbolic groups. To every hyperbolic group $G$ one can associate a boundary at infinity $\dee_\infty G$, a compact Hausdorff topological space equipped with a natural class of visual metrics. 
If $H$ is a finite index subgroup of a group $G$, then $G$ is hyperbolic if and only if $H$ is hyperbolic. Moreover, in this case $\dee_\infty G=\dee_\infty H$, see, e.g., \cite[Section~2]{KB01}.

Every element $g$ of a hyperbolic group $G$ induces a quasisymmetric self-map $\hat g$ of the boundary $\dee_\infty G$. This is a special case of~\cite[Theorem~6.5]{BS00};  see Subsection~\ref{SS:Qr} below for the definition of quasisymmetric maps. 
If $G$ is non-elementary, i.e., not finite or virtually cyclic, the boundary $\dee_\infty G$ is perfect, and the action of $G$ on $\dee_\infty G$ has finite kernel, called \emph{ineffective kernel}~\cite[Ch.~8, 36.-Cor.]{gh}.

Every convex-cocompact Kleinian group is hyperbolic, see, e.g.,~\cite{BH99}. 
A hyperbolic group $G$ acts as a uniform convergence group on its boundary $\dee_\infty G$~\cite[Proposition~1.13]{BB96}. Cocompactness on triples in the context of hyperbolic groups is equivalent to  existence of $\epsilon>0$ such that if $o, p$, and $q$ are arbitrary distinct points in $\dee_\infty G$, then there exists $g\in G$ so that 
$$
d(g(o),g(p)),\ d(g(o),g(q)),\ d(g(p),g(q))\ge\epsilon,
$$
where $d$ is a visual metric.
Conversely, \cite[Theorem~0.1]{bB98} states that if $X$ is a perfect, metrizable, compact, Hausdorff topological space, and a group $G$ acts on $X$ by homeomorphisms as a uniform convergence group, then $G$ is hyperbolic. Moreover, there is a $G$-equivariant homeomorphism of $X$ onto $\dee_\infty G$. 


%

If $G$ is a hyperbolic group, any element of infinite order in $G$ acts as a \emph{loxodromic} isometry on the Cayley graph of $G$. I.e., there are exactly two points $g^+$ and $g^-$ in $\dee_\infty G$ fixed by $g$, that are given by $g^+=\lim_{n\to\infty} g^n$ and $g^-=\lim_{n\to\infty}g^{-n}$. The points $g^\pm$ are called \emph{poles} corresponding to the infinite cyclic group $\left<g\right>$.
The set of poles of loxodromic elements of a hyperbolic group $G$ is dense in $\dee_\infty G$; see, e.g., \cite[Proposition~4.2]{KB01}. 


\subsection{Quasiregular and related maps}\label{SS:Qr}
Recall that a non-constant orientation preserving  continuous map $f =(f_1,f_2)\: U \to \hC$
defined on an open set $U\subseteq \hC$ is called \emph{quasiregular} if $f$ is in the Sobolev space $W_{\rm loc}^{1,2}$
and there exists a positive constant $K$ such that in local coordinates the formal differential matrix $Df =(\partial f_j/\partial x_i)$
satisfies
$$
||Df(z)||^2\leq K\det(Df)(z)
$$
for almost every $z \in U$.  The assumption that $f\in W_{\rm loc}^{1,2}$ means that the first distributional partial derivatives of $f$ are locally in the Lebesgue space $L^2$. 
If $f$ is assumed to be a homeomorphism, it is called a \emph{quasiconformal} map.
A result of Reshetnyak states that any quasiregular map $f\: U\to V$ is a branched covering. This means that $f$  is an open map and for each $w\in V$ the preimage is 
a discrete subset of $U$. The \emph{critical set} of $f$, denoted $\crit(f)$, is the set of all points in $U$ near which $f$ is not a local homeomorphism. The set $\crit(f)$ is necessarily discrete. 

We say that a map $f\: A\to B$ from an arbitrary set $A\subseteq \hC$ to $B\subseteq\hC$ is \emph{quasiregular} if it is open (in relative topology) and is the restriction to $A$ of a quasiregular map $F\: U\to \hC$ defined on an open set $U\subseteq \hC$ that contains $A$. We adapt the same terminology for quasiconformal maps. 

A homeomorphism $f\: X\to Y$ between metric spaces $(X, d_X)$ and $(Y, d_Y)$ is called \emph{quasisymmetric} if there exists a homeomorphism 
$$
\eta\:[0,\infty)\to[0,\infty)
$$ 
such that for any triple of distinct points $o, p$, and $q$, we have
$$
\frac{d_Y(f(p), f(o))}{d_Y(f(q),f(o))}\le\eta\left(\frac{d_X(p,o)}{d_X(q,o)}\right).
$$
The homeomorphism $\eta$ is called a \emph{distortion function} of $f$, and if we want to emphasize it, we say that $f$ is $\eta$-\emph{quasisymmetric}. 

Every quasisymmetric map between domains in $\hC$ is quasiconformal.  The Egg yolk principle~\cite[Theorem~11.14]{jH01} gives a partial converse. Let $B(p,r)$ stand for a disc in $\C$ centered at $p$ of radius $r$. If $f\: B(p, 2r)\to\C$ is $K$-quasiconformal, then the restriction $f|_{B(p, r)}$ is $\eta$-quasisymmetric onto its image with $\eta$ that depends only on $K$.

%

\subsection{Statement of results}\label{SS:SR} 
If $X$ is an arbitrary metric space, in what follows we denote by $B(p,r)$ an open ball in $X$ centered at $p$ of radius $r>0$.
A \emph{Schottky set} $S$ is a compact subset of $\hC$ whose complement is a union of at least three open geometric discs whose closures do not intersect. The boundary circles of the complementary discs are called \emph{peripheral circles}.
If a Schottky set $S$ has empty interior, as is typical in what follows, it is homeomorphic to the standard Sierpi\'nski carpet. 

The main result of this paper is the following theorem.

\begin{theorem}\label{T:Main}
Suppose that $G$ and $\tilde G$ are Kleinian groups whose limit sets $S$ and $\tilde S$, respectively, are  Schottky sets. We assume that the actions of $G$ on $S$ and $\tilde G$ on $\tilde S$ are  cocompact on triples. If $f\: A\to \tilde S$ is a quasiregular map defined on an open (in relative topology) connected subset $A$ of $S$, then $f$ has to be the restriction of a M\"obius transformation that takes $S$ onto $\tilde S$. In particular, $f$ is injective.
\end{theorem}

The fundamental groups of compact hyperbolic 3-manifolds with non-empty totally geodesic boundaries satisfy the assumptions of Theorem~\ref{T:Main}. 

Let $G$ be a hyperbolic group whose boundary at infinity $\dee_\infty G$ is a  Sierpi\'nski carpet. If $\dee_\infty G$  embeds quasisymmetrically into $\hC$, then, according to~\cite[Corollary~1.2 combined with Proposition~1.4]{mB11}, there exists a quasisymmetric map 
$$
\beta\:\dee_\infty G\to S, 
$$
where $S$ is a Schottky set. The peripheral circles of the Schottky set $S$ would necessarily \emph{occur on all locations and scales}~\cite[Proposition~1.4]{mB11}. This means that there exists a constant $C>0$ such that for every $p\in S$ and every $0<r\le2$ there exists a peripheral circle $J$ of $S$ that intersects $B(p,r)$ and such that
$$
\frac{r}{C}\le r_J\le Cr,
$$ 
where $r_J$ is the radius of $J$. The constant 2 in the above definition is the diameter of the standard 2-sphere. This notion is invariant under quasisymmetric maps and it is stronger than local porosity stated in Section~\ref{S:RelS}.
Thus, \cite[Theorem~1.2]{BKM07} implies that the action of $G/K$ on $\dee_\infty G$, where $K$ is the ineffective kernel, is conjugated by $\beta$ to an action of a Kleinian group $G'$ on $S$. The Schottky set $S$ is necessarily the limit set of $G'$. If we assume that $G$ is torsion-free, then $K$ is trivial, and hence the group $G$ acts on $S$ cocompactly on triples. Therefore we have the following consequence of Theorem~\ref{T:Main}.    

\begin{theorem}\label{T:Liouv}
Suppose that $G$ is a torsion-free hyperbolic group whose boundary at infinity $\dee_\infty G$ is a Sierpi\'nski carpet. We assume that $\dee_\infty G$ endowed with a visual metric embeds quasisymmetrically into the standard 2-sphere. Then there exists a group $H$ that contains $G$ as a finite index subgroup and that has the following property. Given any two open connected subsets $A$ and $B$ of $\dee_\infty G$ and a quasisymmetric map $f\: A\to B$, there exists $h\in H$ such that $\hat h|_{A}=f$.    
\end{theorem}

The assumption that the visual boundary embeds quasisymmetrically into the standard 2-sphere is conjecturally true for any hyperbolic group with a Sierpi\'nski carpet boundary~\cite{KK00}. 

The following two corollaries are consequences of Theorems~\ref{T:Main} and \ref{T:Liouv}.

\begin{corollary}\label{C:BoundaryJul}
Assume that  $G$ is a torsion-free hyperbolic group whose boundary $\dee_\infty G$ is a Sierpi\'nski carpet, and let $f$ be an arbitrary rational map. Then $\dee_\infty G$ and the Julia set $\Ju(f)$  of $f$ are not quasisymmetrically equivalent. 
\end{corollary}


\begin{corollary}\label{C:SelfSim}
Suppose that $C$ is a Sierpi\'nski carpet embedded in the standard 2-sphere and that has the following property. There exist an open subset $A$ in $C$ that is not dense in $C$, and a quasisymmetric map $f\: \bar A\to C$, where $\bar A$ denotes the closure of $A$. Then $C$ cannot be quasisymmetrically equivalent to the boundary at infinity of any torsion-free hyperbolic group $G$.
\end{corollary}

The following is a family of Sierpi\'nski carpets that satisfy the assumptions of Corollary~\ref{C:SelfSim}. Compare this to the discussion immediately preceding Theorem~1.3 in~\cite{BM13}.

\medskip
\no
{\bf Example.}
By a \emph{self-similar} carpet $C$ we mean a Sierpi\'nski carpet obtained in the following way. Start with a square in the plane and consider its tiling by finitely many subsquares. Next, remove the interiors of some of the subsquares in the tiling making sure that the following two conditions are satisfied: the closures of any two subsquares whose interiors are removed do not intersect; the closure of any subsquare whose interior is removed is disjoint from the boundary of the original square.  Further,  perform the same operations on the subsquares that remain, i.e., consider the same tiling of each of these subsquares and remove the interiors of the squares of each such tiling so that the combinatorics is the same as in the previous step. Continue the process indefinitely. 

For a quasisymmetric map $f$ in Corollary~\ref{C:SelfSim} we can choose the affine rescaling map from any of the non-trivial scaled copies of $C$ onto $C$. 


\medskip\noindent 
\textbf{Acknowledgments.} The author thanks Mario Bonk and Ilya Ka\-po\-vich for many useful conversations. 

\section{Relative Schottky sets and Schottky maps}\label{S:RelS}
\no
A \emph{relative Schottky set} $S$ in a domain $D\subseteq\hC$ is the residual set obtained by removing from $D$ open geometric discs whose closures are contained in $D$ and are pairwise disjoint. More precisely, we assume that there exists an index set $I$ that consists of at least three elements, and such that 
$$
S=D\setminus\cup_{i\in I} B_i,
$$
where $B_i,\ i\in I$, are open geometric discs with closures $\bar B_i,\ i\in I$, contained in $D$, and $\bar B_i\cap \bar B_j=\emptyset,\ i\neq j$. If $D=\hC$, we recover Schottky sets.

\begin{lemma}\label{L:Jord}
Let $S$ be a Schottky set and $p$ be a point in $S$ that does not belong to any of the peripheral circles of $S$. Then for every open set $U\subseteq\hC$ that contains $p$, there exists a Jordan curve $C$ with the following properties. The curve $C$ is contained in $S$, it does not intersect any of the peripheral circles of $S$, the Jordan domain $D\subseteq\C$ bounded by $C$ contains $p$ and its closure $\bar D$ is contained in $U$.  In particular, $S'=S\cap D$ is a relative Schottky set.
\end{lemma}
\no
\emph{Proof.}
This is a simple application of Moore's theorem~\cite{rM25}. The elements of the decomposition space are points as well as the closures of all the complementary discs of the Schottky set $S$. The decomposition space, i.e., the projection of $\hC$ under the map that identifies points that belong to the same complementary component of $S$, is homeomorphic to the sphere. So we can identify it with the standard 2-sphere. 

If $\tilde p$ is the point of the decomposition space that corresponds to $p$, let small $r>0$ be chosen so that the circle centered at $\tilde p$ of radius $r$ does not contain any of the points of the decomposition space that correspond to the complementary components. This is possible because the number of such components is countable.  Moreover, since $p$ does not belong to any of the peripheral circles of $S$, we can choose $r$ so small that the closure of the disc $B(\tilde p, r)$ is contained in the projection of $U$. Now, the Jordan curve $C$ is the preimage of the boundary circle of $B(\tilde p, r)$ under the projection. The stated properties of $C$ are immediate.
\qed

\medskip

Let $S$ and $\tilde S$ be relative Schottky sets and $f\:A\to \tilde S$ be a local homeomorphism defined on an open subset $A$ of $S$. Then $f$ is called a \emph{Schottky map} if it is \emph{conformal}~\cite{sM13}. This means that for every $p\in A$, the derivative of $f$ defined as  
$$
f'(p)=\lim_{q\in A,\, q\to p}\frac{f(q)-f(p)}{q-p},
$$
exists, non-zero, and continuous.
The following lemma, as well as Lemma~\ref{L:QcExt} below,  are consequences  of Lemma~\ref{L:Jord}.

\begin{lemma}\label{L:QrS}
Suppose that $S$ and $\tilde S$ are Schottky sets of measure zero, and let $f\: A\to\tilde S$ be a quasiregular map defined on an open set $A\subseteq S$. Then $f$ restricted to $A\setminus\crit(f)$ is a Schottky map.
\end{lemma}
\no
\emph{Proof.}
It is enough to show that $f$ restricted to some neighborhood of every point of $A\setminus\crit(f)$ is a Schottky map. 

Assume first that $p\in A\setminus\crit(f)$ does not belong to any of the peripheral circles of $S$. Let $r>0$ be chosen so small that $f$ can be extended to a quasiconformal map $F$ in the disc $B(p, 2r)$. The Egg yolk principle implies that $F$ is quasisymmetric in $B(p,r)$. Lemma~\ref{L:Jord} gives a Jordan domain $D$ that contains $p$, whose closure is contained in $B(p, r)$, and whose boundary $\dee D$ is contained in $S$ and does not intersect any of the peripheral circles of $S$. 

Let $S'=S\cap D$, a relative Schottky set in the Jordan domain $D$. The set $S'$ has measure zero since $S$ does. The restriction of $f$ to $S'$ is a quasisymmetric map $f_D$ from $S'$ onto the relative Schottky set $\tilde S'=\tilde S\cap F(D)$ in the Jordan domain $F(D)$.  Now, \cite[Theorem~1.2]{sM10} implies that $f_D$ is a Schottky map. 

The case when $p\in A\setminus\crit(f)$ does belong to a peripheral circle of $S$ can be reduced to the previous one as follows. Since $f$ is open, it sends every peripheral circle of $S$ intersected with $A$ to a peripheral circle of $\tilde S$. Thus, if $p$ belongs to a peripheral circle $J$ of $S$, we can use the Schwarz reflection principle to extend $f$ across $J$. We denote the peripheral circle of $\tilde S$ that contains $f(p)$ by $\tilde J$. Let us denote the reflection in $J$ by $r$ and the reflection in $\tilde J$ by $\tilde r$. The doubles $S_d=S\cup r(S)$ and $\tilde S\cup\tilde r(\tilde S)$ of $S$ and $\tilde S$, respectively, are clearly Schottky sets. The extension of $f$ to $A\cup r(A)$, defined on $r(A)$ as $\tilde r\circ f\circ r$, is still quasiregular since  geometric circles are removable for quasiconformal maps. But $p$ does not belong to any of the peripheral circles of $S_d$, and so we are in the previous case.
\qed

\medskip

A relative Schottky set $S=D\setminus\cup_{i\in I} B_i$ is said to be \emph{locally porous at a point} $p\in S$ if there exist a neighborhood $U$ of $p$, a positive constant $r_0$
and a constant $C\ge1$, such that for every $q \in S \cap U$ and each $r$ with
$0 < r \le r_0$, there exists $i\in I$ with $B(q, r) \cap B_i \neq\emptyset$ and
$$
r/C \le r_i \le Cr,
$$
where $r_i$ the the radius of $B_i$. A relative
Schottky set $S$ is called \emph{locally porous} if it is locally porous at every point. Local porosity is implied by the property of peripheral circles to occur on all locations and scales. 

The following two theorems are proved in~\cite{sM13} and will be used in what follows.

\begin{theorem}\cite[Corollary~4.2]{sM13}\label{T:Un}
Let $S$ be a locally porous relative Schottky set in $D\subseteq \C$,  and suppose that $A$ is an open connected subset of $S$. 
Let $f$ and $g\: A\to \tilde S$ be Schottky maps into a relative Schottky set $\tilde S$ in a domain $\tilde D$, and consider the set 
$$
E=\{p\in S\cap U\: f(p)=g(p)\}.
$$  
Then $E=A$, provided $E$ has an accumulation point in $A$. 
\end{theorem}


\begin{theorem}\label{T:Der}\cite[Theorem~5.2]{sM13}
Let $S$ be a locally porous relative Schottky set in a  domain $D\subseteq\C$, and $p\in S$ be an arbitrary point. 
Suppose that $U\subseteq D$ is an open neighborhood of $p$ such that $S\cap U$ is connected, and assume that there exists a Schottky map 
$f\: S\cap U\to S$ with $f(p)=p$ and $f'(p)\neq1$.
Let $\tilde S$ be a relative Schottky set in a domain $\tilde D$, and let $(h_k)_{k\in\N}$ be a sequence of Schottky maps $h_k\: S\cap U\to \tilde S$.  We assume that for each $k\in\N$, there exists an open set $\tilde U_k$ so that the map $h_k\: S\cap U\to \tilde S\cap \tilde U_k$ is a homeomorphism, and the sequence $(h_k)$
converges locally uniformly to a  homeomorphism $h\:S\cap U\to \tilde S\cap\tilde U$, where $\tilde U$ is an open set. 
Then there exists $N\in \N$ such that $h_k=h$ in $S\cap U$ for all $k\geq N$.
\end{theorem}

\section{Proof of Theorem~\ref{T:Main}}
\no 
As mentioned in Subsection~\ref{SS:HG}, according to~\cite[Theorem~0.1]{bB98}, the groups $G$ and $\tilde G$ are hyperbolic. Moreover, the limit sets $S$ and $\tilde S$ can be identified with the boundaries at infinity $\dee_\infty G$ and $\dee_\infty\tilde G$, respectively. In particular, \cite[Proposition~1.4]{mB11} gives that the peripheral circles of $S$ and $\tilde S$ appear on all locations and scales, and thus $S$ and $\tilde S$ have measure zero. Lemma~\ref{L:QrS} therefore gives that $f$ restricted to $A\setminus \crit(f)$ is a Schottky map. 
%
%
%

Let $p\in A\setminus\crit(f)$ be an arbitrary point that does not belong to any of the peripheral circles of $S$. We will first show that there exists an open neighborhood $D\subseteq\C$ of $p$ and a M\"obius transformation $m$ that takes $S$ onto $\tilde S$, such that $S\cap D\subseteq A$, and the restrictions of $f$ and $m$ to $S\cap D$ coincide. 

As in the proof of Lemma~\ref{L:QrS}, let $r>0$ be chosen so that $f$ extends to a quasiconformal map $F$ in $B(p,2r)$.
Lemma~\ref{L:Jord} gives a Jordan domain $D$ such that $p\in D$, the boundary $\dee D$ is contained in $S$ and does not intersect any of the peripheral circles of $S$, and the closure $\bar D$ is contained in $B(p,r)$.  Let $S'=S\cap D$, a relative Schottky set in $D$.

As discussed in Subsection~\ref{SS:HG}, the set of poles of loxodromic elements of $G$ is dense in $S$. 
Therefore we may choose a loxodromic element $g$ in $G$ with the pole $g^+$ in $S'$. 
By choosing the radius $r$ above sufficiently small and by possibly passing to an appropriate power of $g$, we may assume that $g(D)\subseteq D$.  Let $o, p, q$ be three distinct points in $S'$. For each $n\in\N$, we can use cocompactness of the action of $\tilde G$ on triples to find $\tilde g_n\in \tilde G$ such that $\tilde g_n\circ f\circ g^n(o),\ \tilde g_n\circ f\circ g^n(p)$, and  $\tilde g_n\circ f\circ g^n(q)$ are $\epsilon$-separated for some $\epsilon>0$ independent of $n$. 
Note that each map $\tilde g_n\circ f\circ g^n$ is a Schottky map defined on $S'$ and that has a $K$-quasiconformal extension $F_n$ to $D$, where $K$ does not depend on $n$.  
Standard compactness arguments imply that the sequence $(F_n)$ subconverges, i.e., some subsequence $(F_{n_k})$ converges locally uniformly in $D$ to a homeomorphism $H\: D\to \tilde S\cap \tilde D$ for some open set $\tilde D$. 

Theorem~\ref{T:Der} now gives that the sequence $(\tilde g_{n_k}\circ f\circ g^{n_k})$ stabilizes, i.e., there exists $N\in\N$ such that
$$
\tilde g_{n_k}\circ f\circ g^{n_k}=\tilde g_{n_{k+1}}\circ f\circ g^{n_{k+1}}
$$
in $S'$ for all $k\ge N$.
This gives, in particular, that the equation
\begin{equation}\label{E:Main}
f=\tilde g \circ f\circ g^{n}
\end{equation}
holds in $g^{n}(S')\subseteq S'$ for some $\tilde g\in\tilde G$ and $n\in\N$. 

Equation~\eqref{E:Main} allows us to inductively extend $f$ to $S\setminus\{g^-\}$. Indeed, $\{g^{-nk}(D)\}_{k\in\N}$ forms an increasing sequence of domains whose union is $\hC\setminus\{g^-\}$. Moreover, 
$$
g^{nk}(S')\subseteq g^{n(k-1)}(S')\subseteq\dots\subseteq S'\quad {\rm and}\quad \cup_{k\in\N}g^{-nk}(S')=S\setminus\{g^-\}.
$$
Since $f$ extends to a quasiconformal map in $D$ and $\dee D\subseteq S$,   Equation~\eqref{E:Main} also gives a way to extend $f$ quasiconformally to $\hC\setminus\{g^-\}$. Every orbit $\mathcal O$ of  complementary  components of $S$ under the action of the cyclic group generated by $\langle g^n\rangle$ has an element contained in $D$. For every orbit $\mathcal O$ we choose one such element, denoted $o$,  and extend $f$ into the interior of $o$ as the restriction of the quasiconformal extension $F$ of $f$. Using   Equation~\eqref{E:Main}, we can now extend $f$ into every complementary component of $S$. The extended map is quasiconformal in $\hC\setminus\{g^-\}$. This is essentially~\cite[Lemma~9.1]{sM10}; it is not crucial that the Ahlfors--Beurling extension into every complementary component was used there.
A point is removable for quasiconformal maps, the map $f$ therefore  has a quasiconformal extension $F_g$ to the whole sphere $\hC$. 

Because the standard sphere is a Loewner space (see~\cite{HK98} for the definition), the map $F_g$ is quasisymmetric, and hence its restriction $f_g$ to $S$ is also quasisymmetric.  But $S$ has measure zero, and therefore~\cite[Theorem~1.2]{BKM07} gives that $f_g$ is the restriction of a M\"obius transformation $m$ that takes $S$ onto $\tilde S$. Since $f=f_g$ in $S\cap D$, the restrictions of $f$ and $m$ to $S'=S\cap D$ coincide.

The desired conclusion of Theorem~\ref{T:Main} is now a consequence of Theorem~\ref{T:Un}.
Indeed, as mentioned above, by Lemma~\ref{L:QrS}, the map $f\: A\setminus\crit(f)\to\tilde S$ is a Schottky map. Clearly, $m\: A\setminus\crit(f)\to\tilde S$ is also a Schottky map. Since $A\cap\crit(f)$ is a discrete subset of $A$, the set $A\setminus\crit(f)$ is still open and connected. From the above we know that the set 
$$
E=\{p\in A\: f(p)=m(p)\}
$$ 
contains an open set. Thus Theorem~\ref{T:Un} implies that 
$$
E\setminus\crit(f)=A\setminus\crit(f). 
$$
The continuity of $f$ and $m$ gives that $E=A$.
\qed

\medskip

\section{Proof of Theorem~\ref{T:Liouv}}
\no
We start by proving the following elementary lemma.

\begin{lemma}\label{L:Finindex}
Suppose that $G$ is a hyperbolic group that acts on a Schottky set $S$ by M\"obius transformations. Then there exists a group $H$ of M\"obius transformations that contains $G$ as a finite index subgroup, and such that if $h$ is a M\"obius transformation that leaves $S$ invariant, then $h\in H$.
\end{lemma}
\no
\emph{Proof.}
Assume for contradiction that there exists an infinite sequence $(h_n)_{n\in\N}$ of M\"obius transformations such that $h_1\notin G$ and $h_{n+1}\notin\langle G, h_1, h_2,\dots, h_n\rangle$, for all $n\in\N$. Let $o, p$, and $q$ be three distinct points in $S$. Since $G$ is cocompact on triples, for each $n\in\N$, there exist an element $g_n\in G$ such that $g_n\circ h_n(o), g_n\circ h_n(p)$, and $g_n\circ h_n(q)$ are $\epsilon$-separated for some $\epsilon>0$. Thus there exists a subsequence $(g_{n_k}\circ h_{n_k})$ that converges to a M\"obius transformation $f$. Since the group of all M\"obius transformations that preserve a Schottky set is discrete (such M\"obius maps have to map peripheral circles to peripheral circles), then 
$$
g_{n_k}\circ h_{n_k}=g_{n_{k+1}}\circ h_{n_{k+1}}
$$ 
for all $k$ large enough. This is a contradiction since $h_{n_{k+1}}$ does not belong to $\langle G, h_1, h_2,\dots,h_{n_k}\rangle$.
\qed

\medskip

Next, we need the following lemma, also a consequence of Lemma~\ref{L:Jord}.

\begin{lemma}\label{L:QcExt}
Assume that $S$ and $\tilde S$ are Schottky sets and $f\: A\to \tilde A$ is a quasisymmetric map between  open subsets $A\subseteq S$ and $\tilde A\subseteq\tilde S$. Then $f$ is quasiconformal, i.e., there exists a quasiconformal map $F$ defined on an open set $U\subseteq\hC$ that contains $A$ and such that $F|_{A}=f$.
\end{lemma}
\no
\emph{Proof.}
Assume that $p\in A$ is an arbitrary point that does not belong to any of the peripheral circles of $S$. Since $A$ is open in $S$, there is an open set $V\in \hC$ such that $A=S\cap V$. Applying Lemma~\ref{L:Jord}, we get a Jordan domain $D$ with the following properties. It contains $p$, its closure is contained in $V$, and its boundary $\dee D$ is contained in $S$ and does not intersect any of the peripheral circles of $S$. 

Let $S'=S\cap D$, a relative Schottky set in $D$. The map $f$ takes it to the relative Schottky set $\tilde S\cap\tilde D$, where $\tilde D$ is the Jordan domain bounded by $f(\dee D)$ and that contains $f(p)$. Thus, \cite[Lemma~9.1]{sM10} implies that $f$ has a quasiconformal extension to $D$. We recall that this extension is obtained by using the Ahlfors--Beurling extension into every complementary component of $S'$ in $D$.

Assume now that $p\in A$ is a point on a peripheral circle $C$ of $S$, and let $\tilde C$ denote the peripheral circle of $\tilde S$ that contains $f(p)$. Let $r$ and $\tilde r$ denote reflections in $C$ and $\tilde C$, respectively, as in Lemma~\ref{L:QrS}. As in that lemma, we can extend $f$ to a quasisymmetric map on $A\cup r(A)$ using the formula $\tilde r\circ f\circ r$ on $r(A)$. Since $p$ does not belong to any of the peripheral circles of the Schottky set $S_d=S\cup r(S)$, the above argument allows us to extend $f$ to a quasiconformal map in a neighborhood of $p$.

We conclude that it is possible to extend $f$ to a quasiconfomal map in a neighborhood of every point $p\in A$. Moreover, if a peripheral circle $C$ of $S$ belongs to any two such neighborhoods, we can always make sure that the two extensions of $f$ inside this peripheral circle are the same. Indeed, the Ahlfors--Beurling extension inside $C$ can always be normalized by three points that only depend on $C$ and not a neighborhood. Therefore all the extensions agree on the intersections,  and thus they produce a quasiconformal extension of $f$ on an open set $U\subseteq \hC$ that contains $A$.
\qed

\medskip
\no
\emph{Proof of Theorem~\ref{T:Liouv}.}
%
As discussed in Subsection~\ref{SS:SR}, we may assume that $S=\dee_\infty G$ is a Schottky set, and $G$ is a Kleinian group that acts on its limit set $S$ cocompactly on triples. Lemma~\ref{L:Finindex} then implies that there is a group $H$ of M\"obius transformations that contains $G$ as a finite index subgroup, and such that if $h$ is a M\"obius transformation that leaves $S$ invariant, then $h\in H$.

According to Lemma~\ref{L:QcExt}, every quasisymmetric map $f\: A\to B$ between two open subsets $A$ and $B$ of $S$ is quasiconformal, and hence quasiregular. We can therefore apply Theorem~\ref{T:Main} to conclude that $f$ is the restriction to $A$ of a M\"obius transformation $h$ that preserves $S$. 
Hence  $h\in H$ and the proof of Theorem~\ref{T:Liouv} is complete. 
\qed

\medskip

\section{Proof of Corollaries~\ref{C:BoundaryJul} and~\ref{C:SelfSim}}

\no
\emph{Proof of Corollary~\ref{C:BoundaryJul}.}
We assume for contradiction that $\dee_\infty G$ (endowed with a visual metric) and $\Ju(f)$ are quasisymmetrically equivalent, i.e., there exists a quasisymmetric map $\alpha\: \Ju(f)\to\dee_\infty G$. This implies, in particular, that $\dee_\infty G$ can be quasisymmetrically embedded into the standard 2-sphere, and $\Ju(f)$ is a Sierpi\'nski carpet.

As in Subsection~\ref{SS:SR}, there exist a Schottky set $S$ and a quasisymmetric map $\beta\: \dee_\infty G\to S$.
Since $\beta\circ\alpha$ is a quasisymmetric map and $S$ is a Schottky set, the peripheral circles of $\Ju(f)$ are uniform quasicircles. Thus, \cite[Proposition~5.1]{mB11} gives that the map $\beta\circ\alpha$ has a quasiconformal extension to the whole Riemann sphere. We denote this extension by $\phi$. The conjugate map $g=\phi\circ f\circ\phi^{-1}$ is a quasiregular map that preserves $S$. 

The action of the group $G$ on $S$ is by M\"obius transformations. Thus we can apply Theorem~\ref{T:Main} to arrive at a contradiction. Indeed, according to this theorem the map $g$ must be the restriction of a M\"obius transformation. In particular, it maps each peripheral circle one-to-one onto a peripheral circle. This readily implies that $g$ cannot have any critical points in the discs that are complementary components of $S$. Thus $f$  cannot have critical points either, i.e., it is a M\"obius transformation. This is impossible because the Julia set of $f$ is a Sierpi\'nski carpet.
\qed  

\medskip

\no
\emph{Proof of Corollary~\ref{C:SelfSim}.}
Assume for contradiction that there exist a carpet $C$ satisfying the assumptions of the corollary and a torsion-free hyperbolic group $G$ such that $C$ and $\dee_\infty G$ are quasisymmetrically equivalent. Thus $\dee_\infty G$ embeds quasisymmetrically into the standard 2-sphere and let  $\alpha\: C\to\dee_\infty G$  be a quasisymmetric map. Furthermore,  let $f\:\bar A\to C$ be a quasisymmetric map guaranteed by the assumptions, where $A$ is an open subset of $C$ that is not dense in $C$.  

The  conjugate map $g=\alpha\circ f\circ \alpha^{-1}$ is a quasisymmetric map of $\alpha(\bar A)$ onto $\dee_\infty G$. But $\alpha(A)$, and thus $g(\alpha(A))$, are open subsets of $\dee_\infty G$. Since the restriction of $g$ to $\alpha(A)$ is still a quasisymmetry, Theorem~\ref{T:Liouv} gives that there exists a quasisymmetry $\hat h$ of $\dee_\infty G$ such that $\hat h|_{\alpha(A)}=  g|_{\alpha(A)}$. This leads to a contradiction as follows. Since $g$ and $\hat h$ are continuous on $\alpha(\bar A)$, we get $\hat h|_{\alpha(\bar A)}=  g|_{\alpha(\bar A)}$. However,
$$
g(\alpha(\bar A))=\alpha(f(\bar A))=\alpha(C)=\dee_\infty G,
$$
but $\hat h(\alpha(\bar A))\neq\dee_\infty G$, because $\alpha(\bar A)\neq \dee_\infty G$. This last assertion follows from the fact that $\alpha$ is a homeomorphism and the assumption that $A$ is not dense in $C$.
\qed

\end{document}